\documentclass[11pt]{article}
\usepackage{amssymb}
\usepackage{graphicx}
\usepackage{verbatim}
\font\elevenss=cmss11

\font\eightss=cmss8

\font\sixss=cmss8 at 6pt

\newfam\ssfam
\textfont\ssfam=\elevenss \scriptfont\ssfam=\eightss
 \scriptscriptfont\ssfam=\sixss
\def\ss{\fam\ssfam \elevenss}%
\font\elevenss=cmss11

\font\eightss=cmss8

\font\sixss=cmss8 at 6pt

\newfam\ssfam
\textfont\ssfam=\elevenss \scriptfont\ssfam=\eightss
 \scriptscriptfont\ssfam=\sixss
\def\ss{\fam\ssfam \elevenss}%
\newtheorem {thm}{Theorem}
\newtheorem {lem}[thm]{Lemma}

\newtheorem {cor}[thm]{Corollary}

\def\Cox{\hfill \Box}

\def\xx{{\bf x}}

\def\S{{\cal S}}
\def\E{{\mathbb E}}
\def\tree{{\cal T}}
\def\P{{\mathbb P}}

\def\Z{{\mathbb{Z}}}
\def\F{{\cal{F}}}

\def\|{\, | \, }

\def\Q{{\bf Q}}
\def\ones{\mbox{\ss ones}}
\def\zeros{\mbox{\ss zeros}}
\def\jump{\mbox{\ss jump}}

\def\spb{{\cal S}}
\def\spa{\Xi}
\def\speed{\mbox{\ss spd}}
\def\speedinf{\mbox{\ss inf-spd}}
\def\speedsup{\mbox{\ss sup-spd}}

\begin{document}

\begin{titlepage}

\begin{center}
{\Large \bf The Klee-Minty random edge chain moves with linear speed}
\end{center}
\vspace{5ex}
\begin{flushright}
J\'ozsef Balogh\footnote{Research supported in part by 
National Science Foundation grant \# DMS 0302804}$^,$\footnote{The Ohio
State University, Department of Mathematics,
231 W. 18th Avenue, Columbus, OH 43210, jobal@math.ohio-state.edu},~\\ 
Robin Pemantle \footnote{Research supported in part by
National Science Foundation grant \# DMS 0103635}$^,$\footnote{University
of Pennsylvania Department of Mathematics, 209 S. 33rd Street, Philadelphia, 
PA 19104, pemantle@math.upenn.edu}
\end{flushright}

\noindent{\bf Abstract:}
An infinite sequence of 0's and 1's evolves by flipping each~1
to a~0 exponentially at rate one.  When a~1 flips, all bits to
its right also flip.  Starting from any configuration with finitely
many 1's to the left of the origin, we show that the leftmost~1
moves right with linear speed.  Upper and lower bounds are given on the speed.



\end{titlepage}

\section{Introduction}

\subsection{Motivation}

The simplex algorithm is widely used to solve linear programs.
It works well in practice, though often one cannot prove that it will.
A parameter in the algorithm is the rule that selects one move 
among possible moves (``pivots'') that decrease the objective function.  
Deterministic pivot rules are known to be possibly very far 
from optimal.  For example, consider problems with the number $m$
of constraints of the same order as the dimension, $n$.   
For virtually every deterministic pivot rule there
is a problem for which the algorithm will take exponential time,
although it is conjectured that there exists a descent path
whose length $O(n)$.  Variants of the original argument
by Klee and Minty~\cite{KM72} are cited in~\cite[page 2]{GHZ98}.  

Randomized pivot rules appear to do better.  According to~\cite{GHZ98},
several of the most popular randomized pivot rules appear to have
polynomial -- even quadratic -- running time.  Rigorous and general 
results on these, however, have been hard to come by.  When one restricts 
to a narrow class of test problems, it becomes possible to obtain some
rigorous results.  G\"artner, Hank and Ziegler~\cite{GHZ98} consider
three randomized pivot rules.  Relevant to the present paper are
their results on the {\em random edge} rule, in which the next move
is chosen uniformly among edges leading to decrease the objective function.
They analyze the performance of this rule on a class of linear programs,
the feasible polyhedra for which are called {\em Klee-Minty cubes},
after~\cite{KM72}.  Such cubes are good benchmarks because they cause
some pivot rules to pass through a positive fraction of the vertices.
They prove that the expected run time is quadratic, up to a possible
log factor in the lower bound:
\begin{thm}[GHZ] \label{th:GHZ}
The expected number, $E_n$ of steps taken by the random edge rule,
started at a random vertex of a Klee-Minty cube, is bounded by
$$\frac{n^2}{4 H_{n+1} - 1} \leq E_n \leq {n+1 \choose 2} \, .$$
Here, $H_n = \sum_{j=1}^n \frac{1}{j} \sim \log n$ is the $n^{th}$ 
harmonic number.
\end{thm}
Their lower bound rules out the possibility that $E_n \sim n \,
{\rm polylog}(n)$ which was twice conjectured by previous 
researchers~\cite[page 29]{PS82},~\cite{Kel81}.  They guess 
that the upper bound is the correct order of magnitude, 
and state an improvement in the upper bound from $(1/2) n^2$ 
to $0.27 \ldots n^2$, whose proof is omitted.  

The method of analysis in~\cite{GHZ98} is via a combinatorial model 
due to~\cite{AC78} which describes the progress of the algorithm
as a random walk on an acyclic directed graph.  In their model,
vertices are bijectively mapped to sequences of 0's and 1's of length $n$,
and each move consists of flipping a 1 (chosen uniformly at random)
to a 0, and simultaneously flipping all bits to the right of the 
chosen bit.  It was in this form that the problem came to our
attention.  Indeed, the remainder of the paper is framed in terms
of a variant of this model, which we find to be an intrinsically
interesting model.  Our main result, Theorem~\ref{th:1}, closes 
the gap left open in~\cite{GHZ98}, proving that the upper bound 
is sharp to within a constant factor and obtaining upper and lower 
bounds differing by a factor of less than~2.  We have moved the 
model to continuous time and made $n$ infinite, since from our view 
as probabilists this is the most natural way to frame such a model.  
Nevertheless, our results apply to the setting of~\cite{GHZ98} as 
well (Corollary~\ref{cor:GHZ}).  Although the model seems simple, 
we remark that we were unable to prove many things about the model, 
including whether certain limits exist.

\subsection{Statement of the model}

The one-dimensional integer lattice is decorated with 0's and 1's,
arbitrarily except that there must be some point to the left of which
lie only 0's.  Each site has a clock that goes off at times
distributed as independent mean-one exponentials.  When a clock rings,
if there is a~0 there nothing happens, but if there is~1 there, 
then it and all (infinitely many) of the sites to the right
flip as well -- 0's become 1's and 1's become 0's.  Later
arguments will use random variables involved in the construction of
this continuous time Markov chain, so we give a formal
construction as follows.

Let $\spb$ be the subset of $\{ 0 , 1 \}^\Z$ consisting of sequences
of 0's and 1's that have a leftmost 1 (equivalently, have finitely many 
1's to the left of the origin).  Let $\{ N(j,t) : t \geq 0 \}_{j \in \Z}$
be a collection of IID Poisson counting processes, that is, step
functions increasing by 1 at random times, which we denote $\xi_{j,i}$,
having independent increments distributed as exponentials of mean 1.
As usual, the filtration is defined by letting $\F_t$ denote the 
$\sigma$-field generated by $\{ N(j,s) : s \leq t \}$.  For any
$\omega \in \spb$, we define a Markov chain starting from $\omega$ as
follows.  If $\omega$ is the zero string, the chain remains at $\omega$.
Otherwise, let $i$ be the position of the leftmost 1 in $\omega$.
First, fix $j > i \in \Z$ and define a Markov chain on configurations
on the index set $\{ i, i+1 , \ldots , j \}$ in the finitary
model: each site $s$ attempts to flip at times $\xi_{s,r}$ for
all $r \geq 0$.  The flip is successful if and only if there is
a 1 in position $s$ at time $\xi_{s,r}^-$, in which case sites
$s+1 , \ldots j$ flip as well.  This is a well defined process
because for each $T$ there are only finitely many jump times $\xi_{s,r}$ 
with $i \leq s \leq j$ and $r \leq T$.  It is clear that for $j' > j$
the Markov chain thus defined on configurations on $[i,j']$
projects down to the Markov chain on configurations on $[i,j]$
by ignoring the bits to the right of $j$.  By Kolmogorov's Extension
Theorem there is an inverse limit as $j \to \infty$.  This is a 
Markov chain, which we denote $\{ \omega_t \}$, taking values in 
trajectories $\omega : [0,\infty) \to \spb$, and having jumps at a
countable set of times $\xi_{s,r}$, at which a 1 flips to a 0
and all bits to its right flip as well.  

We are interested in the speed at which the leftmost 1 drifts to 
the right.  Because this chain converges weakly to the zero state,
it is convenient to renormalize by shifting to the leftmost 1.
Consequently, we define the space $\spa$ to be the space of 
those sequences of 0's and 1's indexed by the nonnegative integers 
that begin with a 1.  We define two functionals $\zeros$ and $\ones$ 
on $\spa$ by letting $\ones (\xx) \geq 1$ be the number of leading 1's:
$$\ones (\xx) := \inf \{ j \geq 1 : \xx(j) = 0 \}$$
and letting $\zeros (\xx) \geq 0$ be the number of successive 
0's after the first 1:
$$\zeros (\xx) := -1 + \inf \{ j \geq 1 : \xx (j) = 1 \} \, .$$
Let $\ones_j$ denote the set $\{ \xx : \ones (\xx) = j \}$ which
form a partition of $\spa$, and let $\{ \zeros_j \}$ denote the 
analogous partition with respect to the values of $\zeros$.  

We now define a Markov chain on the space $\spa$ whose law
starting from $\xx \in \spa$ is denoted $\Q_\xx$.  Pick $\omega
\in \spb$ such that the leftmost 1 of $\omega$ is in some
position $i$ and $\omega (i+j) = \xx (j)$ for all $j \geq 0$.
We construct the Markov chain $\{ X_t \}$ on $\spa$
as a function of $\{ M_t \}$ as follows.  First, define
$\sigma_0$ to be 0, and $i_0$ to be the position of the leading 1 
in $\omega$.  Now recursively we let $\sigma_n$ be the first time 
after $\sigma_{n-1}$ for which $N(i_{n-1} , \cdot)$ increases.
For $\sigma_{n-1} \leq t < \sigma_n$ we let $X_t$ be $M_t$
shifted so that $i_{n-1}$ is at the origin (and ignoring negative
indices).  We let $\jump_n$ denote $\ones (\sigma_n^-)$ and we
let $i_n = i_{n-1} + \jump_n$.  

Sometimes, it is convenient to look at the chain sampled at
times $\sigma_n$.  Thus we let 
$$Y_n := X_{\sigma_n}$$
which is now a discrete-time Markov chain.  We let $\P_\xx$
denote the law of this chain starting from $\xx$.

Since the conditional distribution of $\sigma_n - \sigma_{n-1}$ given
$\F_{\sigma_{n-1}}$ is exponential of mean~1, it follows that
$\sigma_n /n \to 1$.  The distance the leading 1 has moved to 
the right by time $t$ is the sum $\sum_{n : \sigma_n \leq t} \jump_n$,
and therefore the average speed $\speed (n)$ up to time $\sigma_n$ 
is the random quantity 
\begin{equation} \label{eq:def speed}
\speed (n) := \sigma_n^{-1} \sum_{j=1}^n \jump_j \sim n^{-1} \sum_{j=1}^n 
   \jump_j \, .
\end{equation}
It is not {\em a priori} clear, nor in fact can we prove that $\speed (n)$ 
has a limit as $n \to \infty$.  Consequently we define
\begin{eqnarray*}
\speedinf & := & {\rm liminf}_n \speed (n) \; ; \\
\speedsup & := & {\rm limsup}_n \speed (n) \; . \\
\end{eqnarray*}
\begin{quote}
{\bf Problem 1:} Show that the limiting speed exists.
\end{quote}
The Markov chain $\{ Y_n \}$ is a time-homogeneous process on a 
compact space, so taking cesaro averages of the marginals, we see 
these must have at least one weak limit.  Any weak limit is a 
stationary distribution.  
\begin{quote}
{\bf Problem 2:} Show that there is a unique stationary distribution
for the chain $\{ Y_n \}$.
\end{quote}
This would imply a positive solution to Problem 1.  In particular, it
would imply that 
$$n^{-1} \sum_{j=1}^n \jump_j  \to \int \jump_1 \, d\pi$$
where $\pi$ is the unique stationary measure; hence the speed
would not only exist but would be almost surely constant independent
of the starting state.  Although we do not have a solution to Problem~2, 
we believe something stronger may hold.
\begin{quote}
{\bf Problem 3a:} Let $T^j \pi$ denote the composition of the
measure $\pi$ with a translation by $j$ bits, e.g., if
$A$ is the event that there is a 1 in position $r$, then 
$T^j \pi (A) = \pi (A^j)$, where $A_j$ is the event that there
is a~1 in position $r+j$.  Prove that $T^j \pi \to M$ where
$\pi$ is the stationary measure and $M$ is IID fair coin flipping.

{\bf Problem 3b:} Prove or disprove that the unique stationary 
measure $\pi$ is equivalent (mutually absolutely continuous) to $M$.
\end{quote}

To see why we believe (3a) to be true, consider a window of
positions of any size $k$, located far to the right of the initial~1.  
Between each time a clock rings in this window there are many times
a~1 turns to a~0 to the left of the window.  Each time this happens,
the bits in the window all flip.  Projecting configuration space onto
what is visible in the window, and again onto a space of half the size
by identifying each configuration with its complement, it seems reasonable 
to approximate the projection by a Markov chain.  Specifically, from 
a state $\{ x , x^c \}$, for each position $i \in [1, \ldots ,k]$, exactly 
one of $x$ or $x^c$ will be able to flip the bit in position $i$, 
so one may imagine the pair $\{ x , x^c \}$ as flipping at rate $1/2$
in every position.  This chain has the uniform distribution as the
unique stationary distribution.

\section{Statement of main result and lemmas}

In this section we state the results that we do know how to prove,
namely bounds on the lim inf and lim sup speeds.  The following
theorem is to be interpreted as referring to the lim inf and lim
sup speeds (until we have proved the speed exists).

\begin{thm} \label{th:1}
The speed of the drift of the leftmost~1 satisfies
$$1.646 < \speedinf < \speedsup < 2.92 \, .$$
\end{thm}

Relating back to the performance of the random edge rule on 
Klee-Minty cubes, we have:

\begin{cor} \label{cor:GHZ}
For sufficiently large $n$, starting from a uniform random 
vertex of the Klee-Minty cube,
$$0.086 n^2 \leq E_n \leq 0.152 n^2 \, .$$
\end{cor}

\noindent{\sc Proof:} G\"artner, Henk and Ziegler consider another
way of counting steps, where instead of choosing an edge at random
among all those decreasing the objective function, they choose an
edge at random from among all edges, but suppress the move if the
edge increases the objective function.  For a vector $\xx$ of~0's 
and~1's of length $n$, let $L^* (\xx) = L^* (\xx , n)$ denote the
expectation of the number $N^{(r)}$ of moves starting from $\xx^r$, 
including the suppressed moves, before the minimum is reached.  
They prove the following identity~\cite[Lemma~4]{GHZ98}.
\begin{equation} \label{eq:GHZ}
E_n = \frac{1}{2n} \sum_{r=1}^n L^* (\xx^r , n)
\end{equation}
where $\xx^r$ is the vector of length $r$ consisting of all~1's.  
Including suppressed moves in the count corresponds in our infinite, 
continuous-time model to counting the number of clock events 
(only among the first $r$ vertices).  Let $T^{(r)}$ be the time it
takes starting from $\xx^r$ to reach the minimum.  By the strong 
law of large numbers, $N^{(r)} / T^{(r)} \sim r$ as $r \to \infty$.  Also
by the strong law, if the liminf and limsup speed are known to be 
in the interval $(a,b)$, then 
$$\frac{r}{b} < T^{(r)} < \frac{r}{a}$$ 
for sufficiently large $r$.  Hence, for sufficiently large $r$,
$$\frac{n r}{b} < N^{(r)} < \frac{n r}{a}$$ 
and plugging into~(\ref{eq:GHZ}) and summing from $r=1$ to $n$ gives
$$\frac{n^2}{4b} < E_n < \frac{n^2}{4a}$$
for sufficiently large $n$.  Plugging in $a = 2.92$ and $b = 1.646$
proves the corollary.   $\Cox$

The lower bound of Theorem~\ref{th:1} is proved in Section~\ref{ss:lower}.  
The lower bound may in principle be improved so as to be arbitrarily near
the actual speed.  For the upper bound, we state some lemmas.  Let 
$$H_k := \sum_{j=1}^k \frac{1}{j}$$ 
denote the $k^{th}$ harmonic number.  
Let $\{ S_n \}$ be a random walk whose increments are equal $k$
with probability $2/((k+1)(k+2))$ for each integer $k \geq 1$.
Let $S := S_{G-1}$ where $G$ is an independent geometric random variable
with mean~2.  
Let $\Theta$ be a random variable satisfying
$$\P (\Theta \geq j) = 1 - F_{\Theta} (j-1) = \sum_{k=1}^\infty 
   \frac{1}{k} \frac{1}{k+j} = \frac{H_j}{j} \, .$$
Assume $\{ S_n \}$ and $\Theta$ are independent of each other and of
$\{ \F_t \}$; denote expectation with respect to $\P_\xx$ by $\E_\xx$ 
and let $\E$ denote expectation with respect to the laws of $S$ 
and $\Theta$.  Analogously with $\zeros (\xx)$ we define the quantity
$$\zeros^* (\xx) := -1 + \inf \{ j > \ones (\xx) : 
   \xx (\ones (\xx)- 1 + j) = 1 \}$$
to be the number of zeros after the first block of ones; 
thus $\zeros^* (\xx) \geq 1$, $\zeros^* (\xx) = \zeros (\xx)$
if and only if $\zeros (x) \geq 1$ and $\zeros (\xx) = 0$ 
if and only if $\ones (\xx) \geq 2$.

\begin{lem} \label{lem:third}
For any $\xx \in \ones_j$,
$$\E_\xx \jump_1 = \sum_{k=1}^j \frac{1}{k} \, .$$
Equivalently, for any $\xx \in \spa$,
$$\E_\xx \jump_1 = H_{\ones (\xx)} \, .$$
\end{lem}

\begin{lem} \label{lem:second}
For any $j \geq 1$, any $\xx \in \zeros_j$, and any integer $L \geq 1$,
$$\P_\xx (\ones (Y_1) \geq L) \leq \P (S + j \geq L) \, .$$
When $\zeros (\xx) = 0$ then 
$$\P_\xx (\ones (Y_1) \geq L) \leq 
   \P (\Theta + \zeros^* (\xx) \cdot B + S \geq L)$$
where $B$ is a Bernoulli with mean $1/2$, and $\Theta, B$
and $S$ are all independent.
\end{lem}

Since $S + \Theta + \zeros^* (\xx)$ is an upper bound for both
quantities $S + \zeros^* (\xx)$ and $S + \Theta + B \zeros^* (\xx)$
appearing as stochastic upper bounds in Lemma~\ref{lem:second}, and
since $H_n$ increases in $n$, we may put this together with 
Lemma~\ref{lem:third} to obtain
\begin{cor} \label{cor:expect 2}
For any $\xx$,
$$\E_\xx \jump_2 \leq \E H_{\Theta + \zeros^* (\xx) + S} \, .$$
$\Cox$
\end{cor}

\begin{lem} \label{lem:eight}
The conditional distribution of $\zeros^* (Y_n)$ given $\F_{\sigma_{n-1}}$
is always bounded above stochastically by the law of $\Theta$.  In other
words, for all $j$,
$$\P (\zeros^* (Y_n) > j \| \F_{\sigma_{n-1}} ) \leq 1 - F_{\Theta} (j)
   \, .$$
\end{lem}

\noindent{\sc Proof of Theorem~\protect{\ref{th:1}} from the lemmas:}
It suffices to show that for any $\xx \in \spa$,
$$\E_\xx \jump_3 \leq 2.92 \, .$$
We simply iterate conditional expectations and compute.  
By the Markov property, and Corollary~\ref{cor:expect 2}, 
\begin{eqnarray*}
\E_\xx \jump_3 & = & \E_\xx \E_{Y_1} \jump_2 \\
& \leq & \E_\xx \left. \left ( \E_y H_{\S + \Theta + \zeros^* (y)} 
   \right ) \right |_{y = Y_1} \, .
\end{eqnarray*}
Since $H_n$ is increasing in $n$ we may use the stochastic upper
bound in Lemma~\ref{lem:eight} for any $\xx$ to see that this is at most
$\E H_{S + \Theta^{(2)}}$ where $\Theta^{(2)}$ is the sum of two
independent copies of $\Theta$.  The upper bound in Theorem~\ref{th:1} 
is completed by computing an upper bound for this.   

The function $H$ is concave and $\Theta^{(2)} \geq 1$, so 
$$H(\Theta^{(2)} + j) - H(\Theta^{(2)}) \leq H(j+1) - 1$$
and we may therefore conclude that 
\begin{equation} \label{eq:linear sum}
\E H_{\Theta^{(2)} + S} \leq \E H_{\Theta^{(2)}} + \E H_{S+1} - 1 \, . 
\end{equation}

To compute the quantity $\E H_{\Theta^{(2)}}$, let $\Theta_1$ and $\Theta_2$ 
independently have the distribution of $\Theta$ and write
\begin{eqnarray*}
\E H_{\Theta_1 + \Theta_2} & = & 1 + \sum_{1 \leq j,k} \frac{1}{j+k} 
   \P (\Theta_1 \geq j) \P (\Theta_2 = k) \\[1ex]
& = & 1 + \sum_{1 \leq j , k} \frac{1}{j+k}
   \frac{H_j}{j} \left ( \frac{H_k}{k} - \frac{H_{k+1}}{k+1} \right ) \\[1ex]
& = & 1 + \sum_{1 \leq j , k} \frac{1}{j+k}
   \frac{H_j}{j} \left ( \frac{H_k}{k} - \frac{H_{k+1}}{k} 
   + \frac{H_{k+1}}{k} - \frac{H_{k+1}}{k+1} \right ) \\[1ex]
& = & 1 + \sum_{1 \leq j , k} \frac{1}{j+k}
   \frac{H_j}{j} \frac{H_{k+1}-1}{k(k+1)}  \, .
\end{eqnarray*}

We may evaluate this numerically as just under~2.

For the last term on the RHS of~(\ref{eq:linear sum}), let 
$\phi (z) = \E z^S = \sum_{n=0}^\infty z^n \P (S = n)$ be 
the generating function for $S$, so $z \phi$ is the generating
function for $S+1$.  For any positive integer $j$ there is an identity
$$\int_0^1 \frac{1 - z^j}{1-z} \, dz = \int_0^1 (1 + \cdots + z^{j-1}) 
   \, dz = H_j \, .$$
Consequently, we may write
\begin{equation} \label{eq:harm}
\E H_{S+1} - 1 = \int_0^1 \frac{z - \E z^{S+1}}{1-z} \, dz = 
   \int_0^1 z \frac{1 - \phi (z)} {1-z} \, dz  \, .
\end{equation}

To compute the generating function $\phi$, first compute the  
generating function $f$ for the increments of $\{ S_n \}$:
$$f(z) = \sum_{k=1}^\infty \frac{2}{(k+1)(k+2)} z^k = \frac{2z - z^2 - 
   2 (1-z) \log \frac{1}{1-z}}{z^2} \, .$$
Then $\phi = 1/(2-f)$ and the integral in~(\ref{eq:harm}) becomes
$$\E H_{S+1} = \int_0^1 \frac{2 z ( \log \frac{1}{1-z} - z)}{3 z^2 - 2z + 
   2 (1-z) \log \frac{1}{1-z}} \, dz \, .$$
One may evaluate this numerically to approximately $0.918797$.
Adding this to the value for $\E H_{\Theta^{(2)}}$ gives a little
under $2.92$.  Rigorous bounds may be obtained with more care.
$\Cox$

\section{Proofs of Lemmas} \label{ss:upper}

\noindent{\sc Proof of Lemma}~\ref{lem:third}:
By definition, each $\xx \in \ones_j$ begins with $j$ 1's in positions
$0 , \ldots , j-1$ followed by a zero.  The evolution of $\{ M_t \}$
decreases the binary representation $\sum_k 2^{-k} \xx (k)$, whence
$M_{\jump_1^-} \in \ones_k$ for some $k \leq j$, that is, there is
always a zero in some position in $[0,j]$.  Furthermore, once there
is a zero in position $k$ for some $k < j$, then there is always
a zero at or to the left of position $k$.  It follows that
$\jump_1$ is equal to the least $k < j$ for which $\xi_{j,1} < \xi_{0,1}$,
that is, for which the clock in position $k$ goes off before the clock
at position 0.  The minimum is taken to be $j$ if there is no such $k$.

It follows that $\P_\xx (\jump_1 = j) = 1/j$ and that for $0 < k < j$,
\begin{equation} \label{eq:distr 1}
\P_\xx (\jump_1 = j) = \frac{1}{k(k+1)} \, .
\end{equation}
To see this, note that $\jump_1 = k$ if and only if $\xi_{0,1}$
is the minimum of the exchangeable variables $\{ \xi_{k,1}
: 0 \leq k < j \}$.  Similarly, for $0 < k < j$, $\jump_1 = k$
if and only if $\xi_{k,1}$ is the minimum of the variables
$\{ \xi_{r,1} : 0 \leq r \leq k \}$ and $\xi_{0,1}$ is the
next least of the values.  Computing expectations via~(\ref{eq:distr 1})
proves the lemma.   $\Cox$

\noindent{\sc Proof of Lemma}~\ref{lem:eight}: Again let us denote
$q = \ones (\xx)$.  We recall that $\P (\jump_1 = k) = 1/q$
for $k=q$ and $1/(k(k+1))$ for $1 \leq k \leq q-1$.  We claim that
for any $k \leq q$,
\begin{equation} \label{eq:cond timers}
\P (\zeros^* (Y_1) \geq j \| \jump_1 = k) \leq \frac{k+1}{k+j} \, .
\end{equation}
If we can show this, then we will have
$$\P (\zeros^* (Y_1) \geq j) \leq \frac{1}{q} \frac{q+1}{q+j}
   + \sum_{k=1}^{q-1} \frac{1}{k(k+1)} \frac{k+1}{k+j} \, .$$
Changing $q$ to $q+1$ increases this by
$$\frac{j-1}{(q+1)(q+j)(q+j+1)}$$
which is nonnegative.  Setting $q = \infty$ then yields the upper 
bound in the lemma, and it remains to show~(\ref{eq:cond timers}).

Observe first that it suffices to show this for $k = q$.  
This is because when $k < q$, the event $\{ \jump_1 = k \}$
necessitates $\xi_{k,1} = \min \{ \xi_{r,1} : 0 \leq r \leq k \}$.
Thus to evaluate $\P (\zeros^* (Y_1) \geq j \| \jump_1 = k)$ we
may wait until time $\xi_{k,1}$, at which point if no bit to
the left of $k$ has flipped yet, the new conditional probability
$\P (\zeros^* (Y_1) \geq j \| \F_{\xi_{k,1}} , \jump_1 = k)$ is 
always at most $(k+1)/(k+j)$ by applying the claim for $q = k$.

Assuming now that $k = q$, we note that the event 
$\{ \jump_1 = k \}$ that we are conditioning on is just the event
$$A := \{ \xi_{0,1} = \min_{0 \leq i \leq k-1} \xi_{i,1} \}$$ 
is the event that the clock at~0 goes off before any clock in
positions~$1 , \ldots , k-1$.  Conditioning on $A$ then makes 
the law of $\xi_{0,1}$ an exponential of mean $1/k$ without 
affecting the joint distribution of $\{ \xi_{r,s} : r > k \}$.  
Now let $m_1$ be the position at time $\xi_{k,1}$ of the 
first~1 to the right of~$k$, and let $t_1$ be the time 
this~1 flips.  Inductively, define $m_{r+1}$ to be the position
of the first~1 to the right of~$k$ after time $t_r$ and let $t_{r+1}$
be the first time after $t_r$ that this~1 flips.  

If the positions $m_r , \ldots , m_r+j-1$ are not filled with ones 
at time $t_{r-1}$ (define $t_0 = 0$) then it is not possible to have
$\zeros^* (Y_1) \geq j$ and $A$ and $t_{r-1} < \xi_{0,1} < t_r$.
That is, one cannot get from fewer than $j$ ones in the first block of
ones to the right of $k$ to at least $j$ ones at the time of the
flip at~0 without having the leftmost one in this block flip.  On 
the other hand, if these $j$ positions are filled with ones at time 
$t_{r-1}$, then 
$$\P (\zeros^* (Y_1) \geq j , \xi_{0,1} < t_r \| A , \F_{t_{r-1}} , 
   \xi_{0,1} > t_{r-1}) \leq \frac{k}{k+j}$$
since the event $\{ \zeros^* (Y_1) \geq j , \xi_{0,1} < t_r \}$ 
requires that the clock at~0 go off before the clocks in 
positions~$m_r , \ldots , m_r + j - 1$ (recall that conditioning
on $A$ has elevated the rate of the clock at~0 to rate~$k$).   
Similarly, $\P (\xi_{0,1} < t_r \| A , \F_{t_{r-1}} , 
\xi_{0,1} > t_{r-1}) = k/(k+1)$.  Therefore,
$$\frac{\P (\zeros^* (Y_1) \geq j , \xi_{0,1} < t_r \| A , \F_{t_{r-1}} , 
   \xi_{0,1} > t_{r-1})}{\P (\xi_{0,1} < t_r \| A , \F_{t_{r-1}}, 
   \xi_{0,1} > t_{r-1})} \leq \frac{k+1}{k+j} \, .$$

By the craps principle, the RHS is then an upper bound for the 
probability of $\zeros^* (Y_1) \geq j$ conditioned only on 
$\jump_1 = k$.   $\Cox$

\noindent{\sc Proof of Lemma}~\ref{lem:second}: Let $\xx \in \zeros_j$
and first assume $j \geq 1$.  We prove the statement for $n=1$, the proof 
for greater $n$ being identical, conditioned on $\F_{\sigma_{n-1}}$. 
It is simple to check whether $\ones (Y_1) = j$.
The bits in positions $1 , \ldots , j$ will remain 0's until the
leading~1 flips at time $\sigma_1$, so the only thing to check is
whether $\sigma_1 = \xi_{0,1}$ is less than or greater than $\xi_{j+1,1}$.
With probability $1/2$, $\xi_{0,1} < \xi_{j+1,1}$ and in exactly this
case $\ones (Y_1) = j$.

Now condition on this inequality going the other way: 
$\xi_{01} > \xi_{j+1,1}$.  Let $t_1 := \xi_{j,1}$. 
Let $j+1 + k_1$ be the position of the first~0 of $\xx$ to the 
right of position $j+1$.  Then at time $t_1$,
the position of the first one to the right of $j+1$ is $j + Z_1$, where
$Z_1$ is the least $l \in [1,k_1-1]$ for which $\xi_{j+1+l,1} < \xi_{j+1,1}$.
If no such $l$ exists, then $Z_1 = k_1$.  We compute $\P_\xx (Z_1 = l)$
as follows.

The variables $\{ \xi_{r,1} : r \in \{ 0 \} \cup [j+1 , j + k_1]
\}$ are independent exponentials.  For $1 \leq l < k_1$, the event
that $\ones (Y_1) \neq j$ and $Z_1=l$ is the intersection of the
event $A$ that $\xi_{j+1,1}$ is less than $\xi_{0,1}$ and
$\xi_{r,1}$ for all $j + 2 \leq r \leq j + l$ with the event
$B$ that $\xi_{j+1+l,1} < \xi_{j+1,1}$. In other words, among $l+2$
independent exponentials, the index of the least and second least
must be $j + 2$ and $j+1$ respectively.  The unconditional
probability of this is $1/((l+1)(l+2))$.  Having conditioned on
the larger event $\{ \xi_{j+1,1} < \xi_{0,1} \}$, the conditional
probability is therefore equal to $2 / ((l+1)(l+2))$.  This holds
for $l < k_1$, where $Z_1 = k_1$ with the complementary
probability.  To sum up, $Z_1$ is distributed as $S_1 \wedge k_1$
where $S_1$ has the distribution of the random walk increments
described in the lemma.

The last step is to invoke the Markov property.  Condition on
$\F_{t_1}$.  The chain from here evolves under the law $\P_{X (t_1)}$.
Iterating the previous argument, there are two cases.  The first case,
which happens with probability $1/2$ is that the clock at the origin
goes off before the next alarm at location $j+1 + Z_1$.  In this case,
$\ones (Y_1) = j + Z_1$.  In the alternative case, we let $t_2$ be
the time at which the clock at location $j+1 + Z_1$ next goes off.
We let $Z_2$ be the number of consecutive 1's at time $t_2^-$
starting from position $j+1 + Z_1$.  Then conditional on $\F_{t_1}$,
$Z_2$ is distributed as $S_1 \wedge k_2$ where $k_2$ is the
number of consecutive 1's at time $t_1$ starting at position
$j+1 + Z_1$.

Iterating in this way, we have the following inductive definitions.
Let $t_0 = 0$.  Let $\tau$ be the least $r$ for which the clock
at the origin goes off after time $t_r$ but before the first alarm
at location $j+1 + \sum_{i=1}^r Z_i$.  For each $r \leq \tau$, we may define
$k_r$ to be the number of consecutive 1's at time $t_{r-1}$
starting at location $j+1 + \sum_{i=1}^{r-1} Z_i$.  We may then define
$t_r$ to be the first time after $t_{r-1}$ that the alarm
at location $j+1 + \sum_{i=1}^{r-1} Z_i$ goes off, and we may define $Z_r$
so that $j+1 + \sum_{i=1}^r Z_i$ is the location of the first zero to
the right of $j+1 + \sum_{i=1}^{r-1} Z_i$ at time $t_r^-$.

The upshot of all of this is that
$$\ones (Y_1) = j + \sum_{i=1}^{\tau} Z_i$$
and that the joint distributions of $\tau$ and $\{ Z_i :
1 \leq i \leq \tau \}$ are easily described.  Conditioned
on $\tau \geq r$ and on $\F_{t_r}$, the probability of $\tau = r+1$
is always $1/2$; as well, $Z_{r+1}$ given $\tau \geq r+1$ and $\F_{t_r}$
is always distributed as a truncation of $S_1$.  We conclude that
$\ones (Y_1)$ is stochastically dominated by the sum of $\tau$
independent copies of $S_1$, hence as $S_{G-1}$.   

Finally, we consider the case where $\zeros^* (\xx) = l > \zeros (\xx) = 0$.
Let $q = \ones (\xx)$, so that $\xx$ begins with $q$ ones, followed by
$l$ zeros, followed by a one in position $q+l$.  A preliminary
observation is that if we begin with a one at the origin, the
position $W(t)$ of the leading one at a later time $t$ is
an increasing function of $t$; hence, if $T_\mu$ is an independent 
exponential with mean $\mu$, the distribution of $W (T_\mu)$
is stochastically increasing in $\mu$.  

Begin by writing
$$\P_\xx (\ones (Y_1) \geq j) = \sum_{k=1}^q \P (\ones (Y_1) \geq j , 
   \jump_1 = k) \, .$$
Let $l^*$ denote the number of zeros consecutively starting from
position $k$ at time $\xi_{k,1}$ if $\jump_1 = k < q$, and 
$l^* = l$ if $k = q$.  In other words, $l^* = \zeros^* (\xx')$
where $\xx'$ is the word at the last time $t$ that $\ones$ 
changes before the leading bit flips ($t = \xi_{\jump_1 , 1}$
if $\jump_1 < q$ and $t=0$ otherwise).  We may then describe 
$\ones (Y_1)$ as $l^* + W$, where $W$ is the number of 
consecutive positions starting from position $\jump_1 + l^*$ 
that turn to zeros between time $t$ and $\xi_{0,1}$.  Now we
break into two cases.

Condition first on $\{ \jump_1 = q \}$.  The time $\xi_{0,1}$ is 
now an exponential of mean $1/q$, and before this time, the bits 
from position $q+l$ onward evolve independently.  We may describe
$\ones (Y_1)$ as $l + W(\xi_{0,1})$, where $W$ is the number of 
consecutive positions starting at $q+l$ which have become zeros
in the time from~0 to $\xi_{0,1}$.  The first part of this
lemma established that when $\xi$ has rate~1, then $W(\xi)
\preceq S$.  Our preliminary observation now shows, 
conditional on $\{ \jump_1 = q \}$, that $W(\xi_{0,1}) \preceq S)$.

Next, let us condition on $\jump_1 = k < q$, obeserving that then 
$l^* \leq q-k$.  In order to have $l^* \geq r$, it is necessary that 
$\xi_{k , 1} = \min \{ \xi_{s , 1} : k \leq s \leq k + r - 1 \}$.  
Having conditioned on $\jump_1 = k$, the distribution of $\xi_{k,1}$ 
becomes an exponential of rate $k+1$, so that the conditional probability
of this clock going off before $r-1$ other conditionally independent
clocks of rate~1 is just $(k+1)/(k+r)$.  Since the event
$\{ \jump_1 = k < q \}$ has probability $1/(k(k+1))$ if $k < q$
and zero otherwise, we may remove the conditioning and sum to get
$$\P (l^* \geq r , \jump_1 < q) \leq \sum_{k=1}^{q-1} 
   \frac{1}{k(k+1)} \frac{k+1}{k+r} \leq \P (\theta^* \geq r) \, .$$
We also still have in this case $W \preceq S$.  

Putting together the cases $\jump_1 = q$ and $\jump_1 < q$, 
we see that $l^* = l$ with probability $1/q$ and otherwise 
$l^* \preceq \Theta$.  The crude bound $1/q \leq 1/2$ gives
$l^* \preceq \Theta + \zeros^* (\xx) \cdot B$.
Since $l^* \in \sigma (\F_t)$ and the bound $W \preceq S$ holds
conditionally on $\F_t$, we arrive at $\ones (Y_1) \preceq
\Theta + \zeros^* (\xx) \cdot B + S$.   $\Cox$

\section{Lower bound} \label{ss:lower}

To arrive at a lower bound, consider a tree whose vertices are 
positive integers, identified with their binary expansions.  The
root is 1, and the children of $x$ are $2x$ and $2x+1$.  Associated
with each node $x$ are a set of $n(x)$ possible transitions, where
$n(x)$ is the number of 1's in the binary expansion of $x$.  The
transitions are to numbers gotten by flipping a 1 and simultaneously
flipping all bits to its right.  Note that all transitions from 
$x$ are to numbers less than $x$.  Let $r(x)$ denote the {\em reward}
if the leading bit of $x$ is flipped, namely one less than the
number of leading 1's in the binary expansion of $x$.  Recursively,
we assign to each $x$ a mean reward and maturation time, $(a(x), b(x))$
as follows.  Fix a set $B$ of {\em bad nodes}, to be determined later.
On first reading, take $B$ to be empty.  For any $x$, let 
$y_1 , \ldots , y_n$ be the possible transitions from $x$.  Let
$(a(x),b(x)) := (0,0)$ if $x \in B$ and otherwise, let
$$(a(x) , b(x)) := \frac{1}{n} \left ( (r(x) , 1) + \sum_{i=1}^n
   (a(y_i) , b(y_i)) \right ) \, .$$
Let $\tree$ be any finite binary rooted subtree, meaning that any
vertex in the subtree has either zero or two children in the subtree.

\begin{lem} \label{lem:reward}
Suppose no leaf of $\tree$ is in $B$.  Then an almost sure lower bound 
for the lim inf speed from any starting configuration is given by 
the minimum over leaves $x$ of $T$ of $1 + a(x)/b(x)$.  
\end{lem}

The lower bound in Theorem~\ref{th:1} will follow from 
Lemma~\ref{lem:reward} together with an implementation of 
the recursion.  Below is some code written in C that implements 
the recursion for a complete binary tree of depth 15, with 
the set B chosen to give a good bound without much trouble.

\pagebreak
{\small \begin{verbatim}
#include <stdio.h>
#define N 18
#define M     524288

main()  {
        float a[M][2];
        int i , j , k, n, min_index, locmin,  power[N+2], bitcount=0, lastcount=0;
	float min_ratio;
       float  locrat;
       float average;

a[0][0] = 0; a[0][1]=0; a[1][0] = 0; a[1][1] = 0; power[0] = 1; power[1]=2; average = 0;
/* the n loop stratifies the recursion by levels of the tree */
for (n=1; n <= N; n++) 
{
  power[n+1] = 2*power[n];
  /* the i loop indexes vertices inside a level */
  for (i=power[n]; i < power[n+1]; i++) 
  {
    bitcount = 0; lastcount=0; a[i][0] = 0; a[i][1] = 0;
    for (j=0; j <= n; j++) 
    {
      k = i^(power[j]); 
      if (k < i) 
      {
        /* k < i iff there is a 1 in position k */
        k = power[j+1]*(i / (power[j+1])) + power[j+1] - 1 - i % power[j+1];
        a[i][0] += a[k][0]; 
        a[i][1] += a[k][1];
        bitcount++;
	}
      else lastcount=bitcount; 
      /* lastcount will eventually say how many 1's after the leftmost zero */
    }
    a[i][1]++;
    a[i][0] += bitcount - lastcount - 1; 
    a[i][0] = a[i][0] / bitcount;
    a[i][1] = a[i][1] / bitcount;
    /* a particular choice of the set B is given in literal form */
    if (  i==2  ||  i==4  ||  i==5  || i==8  || 
          i==9  ||  i==10 ||  i==11 || 
          i==18 || i==19  ||  i==20 ||  i==21 || i==22 || i==23 ||
          i==36 ||  i==37 || i==38 ||
          i==40  ||  i==41 || i==42 || i==43 ||  i==44 || i==45 || i==46 || i==47 ||
          i==73 || i==74|| i==75 ||
          i==80 || i==81 || i==82 || i==83 || i==84 ||  i==85 || i==86 || i==87 ||
          i==88 || i==89  || i==90 || i==91 ||
          i==160 || i==161 || i==162 || i==163 || i==164 || i==165 || i==166 || i==167 ||
	  i==168 || i==169 || i==170 || i==171 || i==172 || i==173 || i==174 ||  i==175 || i==178 ||
	  i==180 || i==181 || i==182 || i==183 ||
          i==324 || i==325 || i==326 || i==328 || i==329 || 
          i==330 || i==331 || i==332 || i==333 || i==334 || i==335 ||
          i==336 || i==337 || i==338 || i==339 || i==340 || i==341 || i==342 || i==343 || 
          i==344 || i==345 || i==346 || i==347 || i==348 || i==349 ||
          i==350 || i==351 ||
          i==361 || i==362 || i==363 ||
          i==656 || i==657 || i==658 || i==659 || i==660 || i==661 ||  i==662 || i==663 ||
	  i==672 || i==673 || i==674 || i==675 || i==676 || i==677 || i==678 || i==679 ||
          i==680 || i==681 || i==682 || i==683 ||
          i==684 || i==685 || i==686 || i==687 || i==688 || i==689 ||
          i==690 || i==691 || i==692 || i==693 || i==694 || i==695 || i==697 ||  i==698 || i==699||
          i==1322 || i==1323 ||
          i==1346 || i==1347 || i==1348 || i==1349 ||
          i==1350 || i==1351 || i==1352 || i==1353 || i==1354|| i==1355 ||
          i==1356 || i==1357 || i==1358 || i==1359 ||
          i==1360 || i==1361 || i==1362 || i==1363 || i==1364 || i==1365 || 
          i==1366 || i==1367 || i==1368 || i==1369 || i==1370 || i==1371 ||
          i==1372 || i==1373 || i==1374 || i==1375 || i==1380 || i==1381 ||
          i==1384 || i==1385 || i==1386 || i==1387 || i==1388 || i==1389 || i==1390 || i==1391 || 
          i==2704 || i==2705 || i==2706 || i==2707 || i==2708 || i==2709 || i==2710 || i==2711 ||
          i==2720 || i==2721 || i==2722 || i==2723 || i==2724 || i==2725 || i==2726 || i==2727 ||
          i==2728 || i==2729 || i==2730 || i==2731 || i==2732 || i==2733 || i==2734 || i==2735 ||
          i==2736 || i==2737 || i==2738 || i==2739 || i==2740 || i==2741 || i==2742 || i==2743 ||
          i==2745 || i==2746 || i==2747 ||
          i==5456 || i==5457 || i==5458 || i==5459 || i==5467 ||
          i==5460 || i==5461 || i==5462 || i==5463 || i==5465 || i==5466  ) 
	  {
		a[i][0] = 0; a[i][1] = 0; 
	  }       
    if(i == (M/2))
      {
	min_index = i;
	min_ratio = a[i][0]/a[i][1];
      }
     
      if(i > (M/2))
      {
	if((a[i][0]/a[i][1]) < min_ratio)
	  {
	    min_index = i;
	    min_ratio = a[i][0]/a[i][1];
	  }
       }
      if  (i > (M/2))  
	   {  
	       average = average + (a[i][0]/a[i][1]);
           }
      if(i == 3) 
      {
	locmin = i;
	locrat = a[i][0]/a[i][1];
      }      
      if(i < (M/2))
	{  if (i > 71)
          {
	    if(a[i][1] > 0)
             { if  ((a[i][0]/a[i][1]) < locrat)
	       {
	        locmin = i;
	        locrat = (a[i][0])/(a[i][1]);
	       }
	     }
           }
        }
  }
  printf("\n");
  }
 printf("Row with minimum ratio is:\n");
 printf("min_index=%i,exp. reward=%f, exp. time=%f, ratio=%f\n",min_index, 
      a[min_index][0], a[min_index][1], a[min_index][0]/a[min_index][1]);
 printf("Row with local minimum ratio is:\n");
 printf("locmin=%i,exp. reward=%f, exp. time=%f, ratio=%f\n", locmin, locrat, average, M);
	}}
\end{verbatim} 
}

A look at the data shows the minimum value of $a(x)/b(x)$ on each level
to be obtained when the binary expansion of $x$ alternates.  In 
particular, the global minimum is at $x = 349525$ and has value $0.646\ldots$,
which proves the lower bound.

\noindent{\sc Proof of Lemma}~\ref{lem:reward}:
Any finite rooted binary subtree induces a prefix rule, that is,
a map $\eta$ from infinite sequences beginning with a 1 to leaves
of $\tree$, defined by $\eta (\xx) = w$ for the unique leaf of 
$\tree$ that is a prefix of $\eta$.  

Given a trajectory of the Markov chain $\{ X_t \}$, define a sequence
of elements of $\tree$ as follows.  Let $x_0 := \eta (X_0)$ be the
prefix of the initial state of the trajectory.  Let $\tau_0 := 0$.
As the definition proceeds, verify inductively that for 
$\tau_k \leq t < \tau_{k+1}$, the string $x_k$ will be an initial
segment of $X_t$.  The recursion is as follows.
Let $\tau_{k+1}$ be the first time after $\tau_k$ that a 1 flips
in the initial segment $x_k$ of $X_t$.  Let $x_{k+1}'$ be the
string gotten from $x_k$ by flipping this bit and all bits to its right.
If $x_{k+1}' \notin B$ and $x_{k+1}'$ is not the zero string then
let $x_{k+1}$ be $x_{k+1}'$, stripped of any leading zeros.  If
$x_{k+1}' \in B$ or $x_{k+1}'$ is the zero string, then let $x_{k+1} = 
\eta (X_{\tau_{k+1}})$.  

Let $\rho_0 , \rho_1 , \ldots$ be the successive times that this
latter transition occurs, that is, successive times $t$ in the 
recursion that $\eta (X_t)$ is computed.  Given $x \in \tree$, 
suppose we begin counting every time the sequence $\{ x_k \}$ 
hits $x$ and stop counting at every time $\rho_k$.  When we are
counting, we count how many times the leading 1 flips, and how
many 1's flip together when this happens, or more precisely, we
keep a cumulative count, each time adding one less than the number 
of 1's that have flipped together.  Let this cumulative count
be denoted $(A(x,t) , B(x,t))$, where $A$ counts flips of the leading 1.
If $x \in B$, then by convention we do not count anything.   

Claim: For each $x \in \tree \setminus B$, if $x$ is visited infinitely 
often then the ratio of $A(x,t)/B(x,t)$ converges as $t \to \infty$
to $a(x) / b(x)$.  Proof: Conditioned on the past, each time we
start counting we are equally likely to transition to each of the
$n(x)$ possible transitions.  The claim then follows from induction
on $x$.  

Each $x$ defines a set of time where $x$ is {\em on},
that is, where we are counting leading~1 flips and rewards for that $x$.
As $x$ ranges over the leaves of $\tree$, the on-timesets for $x$ 
partition $[0,\infty)$.  In each such time set, we have shown 
that the mean number of recorded simultaneous 1-flips per leading 1 flip
is $1 + a(x)/b(x)$ in the limit, as long as the time set is unbounded.
Since the partition is finite, it follows that the lim inf speed
is at least the minimum of the values.   $\Cox$

We remark that the only place this computation is not sharp is when
the number of 1's flipping together exceeds the number recorded, 
because the present knowledge of the prefix was a string of all 1's
and there were more 1's after this that also flipped.  Thus by 
making the tree $\tree$ big enough, even without increasing $B$, 
we can get arbitrarily close to the true value.  

\section{Further observations}

The following argument almost solves Problem~(3a), and perhaps may
be strengthened to a proof.  Lemma~4 of~\cite{GHZ98} is proved by
means of a duality result.  The result is that the probability, 
starting from a uniform random state, of finding a~1 in position 
$r$ after $t$ steps (counting suppressed transitions), is equal to 
half the probability that $\xx^r$ has not reached the minimum yet
after $t$ steps (again counting suppressed transitions).  The
argument that proves this may be generalized by introducing a
simultaneous coupling of the process from all starting states.
The probability, from a uniform starting state, of finding a~1 in
every position in a set $A$ after $t$ steps, is then the expectation
of the function that is zero if the column vector of all 1's is
not in the span of the columns of the matrix whose rows are the
states reached at time $t$ starting at $\xx^r$, as $r$ varies over $A$,
and is $2^{-u}$ if the column vector of all~1's is in the span
and the matrix has rank $u$.  The kernel of the matrix is the
set of starting configurations that reach the minimum by time $t$
(the simultaneous coupling is linear).  Hence, as long as $A$ and $t$ 
are such that the probability of reaching the minimum from any $\xx^r$ 
by time $t$ goes to zero, the rank of the matrix will be $|A|$ and 
the probability of finding all~1's in positions in $A$ at time $t$
will go to $2^{-|A|}$.  In particular, if a window of fixed size
moves rightward faster than the limsup speed, then what one sees
in this window approaches uniformity.  This is not good enough to
imply uniformity of a window a fixed distance to the right of the 
leftmost~1.  

{\bf Acknowledgements:}
We would like to say thanks to Attila P\'or, sharing this problem, 
and to Jiri Matousek giving references about the Klee-Minty cube.

\end{document}